\documentclass[12pt]{amsart}
\pdfoutput=1

\usepackage[mathcal]{euscript}
\usepackage{mathrsfs}
\usepackage{amsfonts,amssymb,amsbsy,amsmath}
\usepackage{graphicx,color}
\usepackage{epsfig}
\usepackage[all,cmtip]{xy}
\usepackage{verbatim}

\usepackage{setspace}

\makeatletter
\let\uppercasenonmath\@gobble
\makeatother

\newcommand{\R}{\mathbb{R}}

\newcommand{\C}{\mathbb{C}}
\newcommand{\Z}{\mathbb{Z}}

\newcommand{\ba}{\begin{array}}
\newcommand{\ea}{\end{array}}
\newcommand{\ra}{\rightarrow}

\theoremstyle{plain}
\newtheorem{theorem}{Theorem}
\newtheorem{lemma}[theorem]{Lemma}
\newtheorem{prop}[theorem]{Proposition}

\theoremstyle{remark}
\newtheorem{remark}{Remark}
\newtheorem{conj}{Conjecture}

\theoremstyle{definition}
\newtheorem*{defi}{Definition}

\begin{document}
\title{Real open books and real contact structures}

\author{Fer\.{\i}t \"{O}ZT\"{U}RK}
\address{Bo\u{g}az\.{\i}\c{c}\.{\i} \"{U}n\.{\i}vers\.{\i}tes\.{\i}, Department of Mathematics, TR-34342
  Bebek, \.Istanbul, Turkey}
\email{ferit.ozturk@boun.edu.tr}

\author{Nerm\.{\i}n SALEPC\.{I}}
\address{Institut Camille Jordan,
Universit\'e Lyon I,
43, Boulevard du 11 Novembre 1918
69622 Villeurbanne Cedex, France}
\email{salepci@math.univ-lyon1.fr}


\begin{abstract}
A real 3-manifold is a smooth 3-manifold together with an orientation preserving smooth involution, called a real structure. 
In this article we study open book decompositions on smooth real 3-manifolds that are compatible with the real structure. We call them real open book decompositions. 
We show that each real open book carries a real contact structure and two real contact structures supported by the same real
open book decomposition are equivariantly isotopic.
We also show that every real contact structure on a closed 3-dimensional
real manifold is supported by a real open book.
Finally, we conjecture that two real open books on a real contact manifold
supporting the same real contact structure are related by 
positive real stabilizations and equivariant isotopy and that the Giroux correspondence applies to real manifolds as well
namely  that there is a one to one correspondence between the real contact 
structures on a real 3-manifold up
to equivariant contact isotopy and the real open books up to positive real stabilization.
Meanwhile, we study some examples of real open books and real Heegaard decompositions in lens spaces. 
\end{abstract}

\maketitle

\section{Introduction and basic definitions}

A {\em real structure} on an oriented $2n$-manifold (respectively $2n-1$-manifold) $X$, possibly with boundary, is defined as
an involution $c_{X}$ on $X$ satisfying the following two conditions: \\
(i) $c_X$  is orientation preserving if $n$ is even and orientation reversing if $n$ is odd; \\
(ii) The fixed point set of $c_X$   is either empty or has  dimension $n$ (respectively $n-1$).

Hence, if $M$ is the oriented boundary of an oriented  $2n$-manifold $X$, a real 
structure on $X$ restricts to a real structure on $M$.
We call a manifold together with a real structure a {\em real manifold} and the fixed point set of the real structure the {\em real part}.
If a contact structure (or a symplectic structure) pulls back to minus itself under a real structure $c$, such structures are called $c$-real
and those manifolds are called $c$-real contact (or $c$-real symplectic respectively) manifolds.
Real algebraic varieties and links of real algebraic isolated singularities supplied with natural structures are examples for real symplectic and contact manifolds.

Every closed oriented 3-manifold admits an open book decomposition \cite{alx} and a positive contact structure \cite{mar}. 
Furthermore the work of E.~Giroux points out that there is a one to one correspondence
between open book decompositions up to positive stabilization and positive contact structures up to isotopy \cite{gi1} (see e.g. \cite{co}, \cite{et} and \cite{vau} for a careful discussion and for the definitions of  open book decompositions and contact structures). 

However not every 3-manifold admits a real structure; in fact the ones which do not admit a real structure are in abundance \cite{pup}.
Nevertheless, once a real structure on a 3-manifold is given we can make appropriate definitions similar to the ones 
above  and show that there are analogous relations. 
For instance, a {\it real open book decomposition} is an open book with the real structure preserving
the page structure and leaving exactly two pages invariant. 
A {\it $c_M$-real (or characteristic)  Heegaard splitting} is a $c_{M}$-invariant Heegaard surface
with the two handlebodies being mapped to each other by $c_{M}$. 
It is known that every real 3-manifold $(M,c_{M})$ admits a characteristic
Heegaard splitting \cite[Proposition~2.4]{nag}.
In Section~\ref{rak} we present examples of real open books and real Heegaard decompositions. 
We state some observations related to real Heegaard decompositions and real open books on lens spaces
in Section~\ref{mercek}. We describe some real lens spaces which cannot have a real Heegaard splitting with genus 1 or 2. 
We also observe that there are (infinitely many)  closed oriented real 3-manifolds 
which cannot have a {\it maximal}  real Heegaard decomposition; in a maximal real Heegaard decomposition,  the real part of the real structure has exactly genus$+1$ connected components (Proposition~\ref{maksimal}).

A {\it positive real  stabilization} of a real open book is performed by attaching a 2-dimensional 1-handle to
the page of an abstract open book and modifying the real structure and monodromy appropriately (see Section~\ref{rstab}).
This boils down to taking a connected sum with a real tight contact $S^3$, which is unique up to equivariant contact isotopy ("Classification of real tight contact structures on solid tori, $S^3$ and $\R P^3$", a preprint of the authors).
With this definition, we conjecture a Giroux correspondence in the existence of a real structure:
there is a one to one correspondence between the $c_{M}$-real  contact structures on a real 3-manifold $(M,c_{M})$ up
to equivariant contact isotopy and the $c_{M}$-real open books on $M$ up to positive real
stabilization and equivariant isotopy (Conjecture~\ref{jiru}). 
In this direction we can show that two real contact structures supported by the same real
open book decomposition are equivariantly isotopic (Proposition~\ref{aynitemas}).
We also show that every real contact structure on a closed 3-dimensional
real manifold is supported by a real open book (Proposition~\ref{tasir}).
We conjecture that two real open books on a real contact manifold
supporting the real contact structure are related by 
positive real stabilizations and equivariant isotopy (Conjecture~\ref{aynikitap}).
Since some details of the proof of the nonequivariant claim is not available to us,
we cannot prove this last statement but believe that it is correct. Note that 
a proof for Conjecture~\ref{aynikitap} will  in turn finish the proof of  Conjecture~\ref{jiru}.

The aim of the present work is to take a step towards using real contact topology in the research on real algebraic surface singularities in the one hand and on involutions on three manifolds in the other. For the former, let us restate that the link manifold of an isolated real algebraic surface
singularity carries naturally a real contact structure, which can be used to make new observations about the singularity. For instance, among other methods, it is possible to use contact topology to prove that the link manifolds of the real singularities $X_n=\{x^{n+1} + y^2 \pm z^2 = 0\}$ are never equivariantly diffeomorphic although they are  stably equivalent for $n$ even.

For the latter path of research above, we have the following programme in mind. In all the examples of real 3-manifolds we know, there is an associated real contact structure. We believe that this holds in general; i.e. every  real 3-manifold has an associated real contact structure. In the proof that we propose, we proceed topologically by constructing a real open book and then conclude by using the main result of the present work.  Now, provided that our belief holds true, it will be possible to obtain restrictions on 3-manifolds for not having orientation preserving involutions. Namely, assuming that there is a real structure $c$ on a given 3-manifold, 
we will take a $c$-real contact structure and use recent tools about contact 3-manifolds
(i.e. Heegaard-Floer homology, the $\Z_2$-action on it, and the contact invariant) to deduce restrictions to the existence of a real structure.

In the sequel, instead of using the term $c$-real, we usually drop the reference to $c$ whenever the real structure is understood.

\textbf{Acknowledgements.} The second author would like to express her gratitude to E.~Giroux who generously share his deep knowledge on the subject; she  would also like to thank J.Y.~Welschinger and  A.~Degtyarev for fruitful discussions. The authors are grateful to the referee for major corrections and helpful suggestions.
The project leading to this article 
is supported by the Scientific and Technological Research Council of Turkey [TUBITAK-ARDEB-109T671].
The second author is supported by the European Community's Seventh Framework Programme ([FP7/2007-2013] [FP7/2007-2011]) under grant agreement no~[258204].

\subsection{Real open books}

\begin{defi} Let $(B, \pi)$ be an open book decomposition of a real 3-manifold $(M, c_{M})${, where $B$ is a closed 1-manifold in $M$ and $\pi:M-B\ra S^1$}. We say that $(B, \pi)$ is a {\em real open book decomposition},  (or shortly, a \emph{real open book}) if $\rho  \circ \pi=\pi \circ c_{M}$ where $\rho: S^1 \to S^1$ is { the reflection on $S^1=\{z\in\C:|z|=1\}$}  induced from the complex conjugation.  
In particular,  we have  $c_{M}(B)=B$. 

An \emph{isomorphism} between two real open book decompositions of $(M,c_{M})$ is a pair of orientation preserving equivariant diffeomorphisms of $(M,c_{M})$ and $(S^1, \rho)$ commuting with the projections.
\end{defi}

{ As orientation preserving reflections on $S^1$ form a single class, fixing $\rho$ as above is not restrictive.} 
By definition, any real open book decomposition of $(M, c_{M})$ has two pages, 
$S_{-}=\pi^{-1}(-1),\, S_{+}=\pi^{-1}(+1)$, which are invariant under the action of the real structure $c_{M}$. The restrictions $c_{-}=c_{M}|_{S_{-}},\, c_{+}=c_{M}|_{S_{+}}$ yield real structures on $S_{-}$ and, respectively, $S_{+}$. These pages together with the inherited real structures are called  \emph{real pages}.  
One of the fundamental properties of the monodromy $f$ of a real open book is that $f=c_{+ }\circ c_{-}$.  Or equivalently, $f^{-1}=c_{-}\circ  f\circ  c_{-}$  as well as $f^{-1}=c_{+}\circ  f\circ  c_{+}$. (Related to the definition above and the argument in the sequel, see \cite{n1}.)
Therefore, a real open book decomposition of $(M, c_{M})$  gives rise to a triple $(S,f, c)$ where $S$ is one of the real pages, $c$ is the inherited real structure on $S$ and $f$ is the monodromy satisfying  $f|_{\partial S}=$id and $f^{-1}=c \circ f \circ c$.
By definition, $c$ acts on $\partial S$, which is the disjoint union of finitely many circles. Hence either $c$ acts on a connected component of $\partial S$
as reflection or $c$ swaps a pair of connected components, reversing the induced orientations
as boundary.

\begin{lemma}\label{yayilanc}
Let $S, f, c$ be as above. Then the real structure $c:S\to S$ extends to a real structure on $M$.
\end{lemma}

\noindent {\it Proof:}
Here we explain how a real structure on a page extends to a real structure on a 3-manifold $M_f=S\times [0,1]/_{\tiny(f(x),0)\sim(x, 1)} \underset{|\partial S|}{\bigcup} S^1\times D^2$, described by the abstract open book $(S,f)$.  As $M_{f}$ is diffeomorphic to $M$, the real structure obtained on $M_{f}$ can be pulled to $M$ by means of a chosen diffeomorphism between $M_{f}$ and $M$.

\begin{figure}[h]
   \begin{center}
    \includegraphics[scale=0.3]{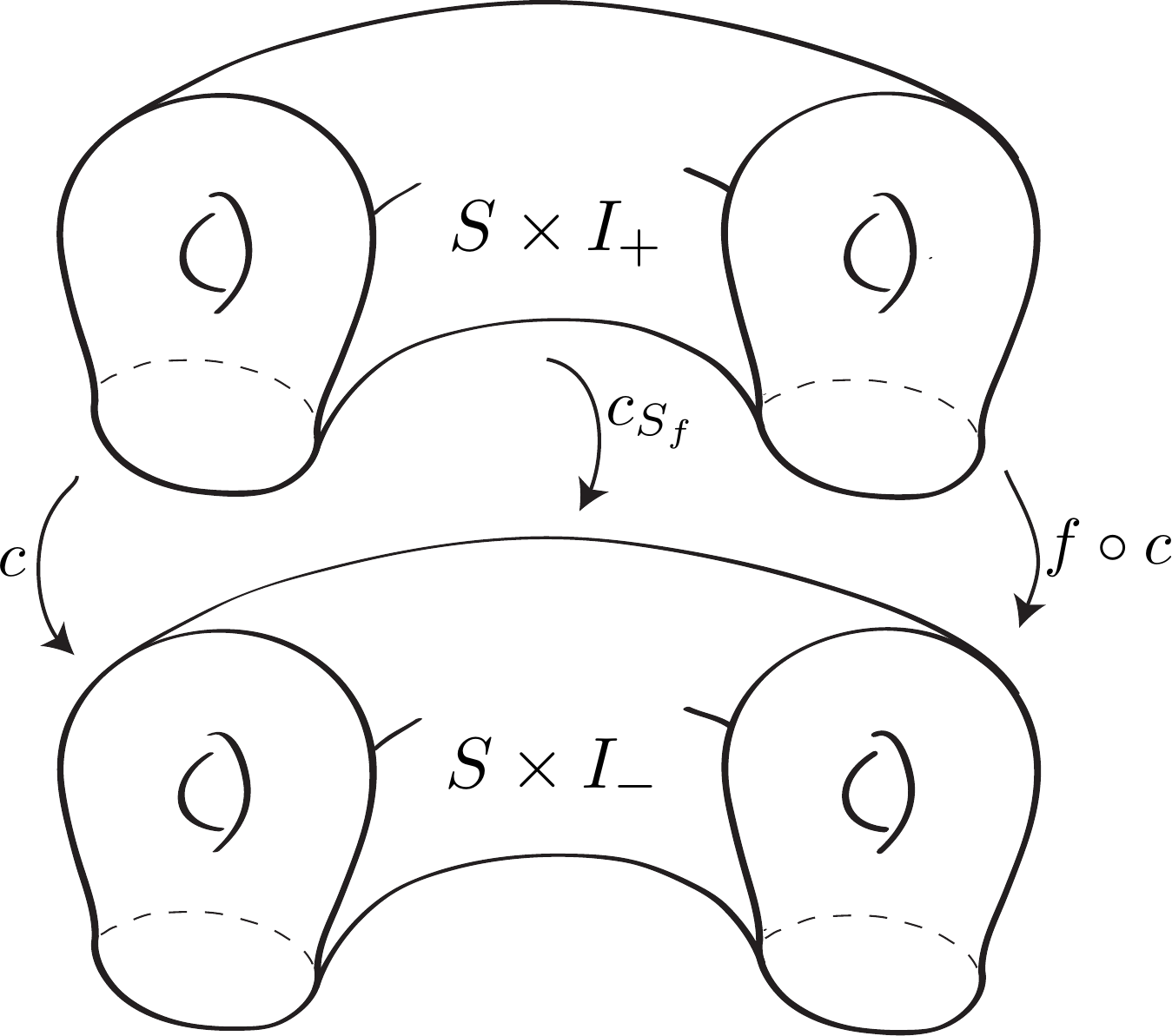}
       \caption{A model for a real open book.}
       \label{boru}
      \end{center}
 \end{figure}

First we describe the real structure on the mapping  torus. 
Instead of the usual description, consider the mapping torus as 
$$S_f=\big((S\times I_+)\cup(S\times I_-)\big)/_{(x,0)\sim (c(x),0) \mbox{ and } (x,1)\sim(f \circ c(x),1)}$$  
for all $x$ in $S$ (see Figure~\ref{boru}). Here $I_\pm$ are copies of $[0,1]$.

Note that the manifold constructed has monodromy $c^{-1}\circ f^{-1} \circ c=f$.

Now consider the map $c_{S_f}:S_f\rightarrow S_f$ which acts as identity 
between the cylinders,
i.e. $(x,t)$ in  $S\times I_+$ is sent to $(x,t)$  in  $S\times I_-$  and vice versa.  After identification $S^1=(-I_+)\cup I_-$ (given by $\theta\mapsto  -t$ if $t\in I_+$ et $\theta\mapsto  t$ if $t\in I_-$) $c_{S_{f}}$ induces a reflection $\theta\mapsto -\theta$ on $S^1$ and hence the projection $S_f \rightarrow S^1$ becomes $(x,\theta)\mapsto  \theta$.

Restricted to each torus component of $\partial S_f$, the map $c_{S_f}$ is  either a rotation by $\pi$ fixing exactly 4 points or it interchanges two tori components.
The various ways to extend  $c_{S_f}$ as a real structure over the 
solid torus neighborhoods of binding components are dictated by the behavior of $c$
on $\partial S$. Either a pair of solid tori are
mapped to each other by $c_{S_f}$ or $c_{S_f}$ extends over a solid torus as rotation by $\pi$.
In each case
the extension is  unique up to isotopy (see e.g. \cite[Lemma~4.4]{hr}) and  the extended map preserves the  page structure.
Therefore, we obtain a real structure on $M_{f}$, described by $(S,f,c)$, preserving the  page structure.
\hfill  $\Box$ \\

Thus, we have the following definition.

\begin{defi} Let $(S,f)$ be an abstract open book, where $S$ is a compact surface with boundary
and $f:S\rightarrow S$ is the monodromy so that $f$ is the identity on $\partial S$.
Suppose that $c$ is a real structure on $S$, i.e. an orientation reversing involution.
An {\em abstract real open book} is a triple $(S,f,c)$ with $f$ satisfying $f \circ c=c \circ f^{-1}$.  

An \emph{isomorphism} of abstract open books is an orientation preserving diffeomorphism of the surface commuting both with the monodromy and with the real structure.
\end{defi} 

Note after the proof of Lemma~1 that whenever we  represent a real manifold with an embedded open book as a real abstract open book, we must allow an equivariant isotopy. 
Note also that once we adopt the model in Figure~\ref{boru} 
so that $c$ is placed on the {\it left}, an (embedded) real open book determines 
an abstract real open book uniquely up to the equivalence above.

\subsection{Positive real stabilization}

\label{rstab}

Recall that a positive stabilization of an abstract open book $(S,f)$ is the abstract open book with page $S'=S\,\cup$ a 1-handle
and with monodromy $\tilde{f} \circ \tau_a$ where $\tilde{f}$ is the trivial extension of $f$ over the 1-handle and $\tau_a$ is a  right-handed Dehn twist along a curve $a$ in $S'$ that intersects the co-core of the 1-handle exactly once.
The type of stabilizations that we need should lead to a real structure on the new abstract open book. Here
is the definition which will provide that. 

\begin{defi}
Let $(S,f,c)$ be an abstract real open book. Let $S'=S\cup H$ where 

\begin{itemize} 
\item either $H$ is a 1-handle with its 
attaching region $c$-invariant,  (in particular if the attaching region is a neighborhood of a pair of real points, then we impose the condition that the real points belong to the same real component);
\item or $H=H_1\cup H_2$ where $H_1$ and $H_2$ are 1-handles with their attaching regions  interchanged
by $c$.
\end{itemize} 
In such cases, $c$ extends uniquely  over $H$ to a real structure, say,
$\tilde{c}$ on $S'$, { up to isotopy through real structures}. Let $\tilde{f}$ denote the extension of $f$ over $S'$ with $\tilde{f}|_H=$id.
We consider, in the former case, a simple closed curve $a$ such that $\tilde{c}(a)=a$ and that $a$ intersects the co-core of $H$ once (existence of 
such an invariant curve $a$ is guaranteed by the imposed condition), while in the latter case, a pair of simple closed curves $a, \tilde{c}(a)$ such that $a$ 
and hence $\tilde{c}(a)$ intersects the co-core of $H_{1}$  and, respectively, of $H_{2}$
once. Depending on the succeeding cases, let $\sigma$ denote either the Dehn twist $\tau_a$ along $a$ or the product $\tau_a \circ  \tau_{\tilde{c}(a)}$. 
Then a {\it positive real stabilization} of the abstract real open book $(S,f,c)$ 
is defined as the real  open book $(S',\tilde{f}\circ \sigma,\tilde{c}\circ \sigma)$. 
\end{defi}

\begin{figure}[h]
   \begin{center}
\includegraphics[scale=0.37]{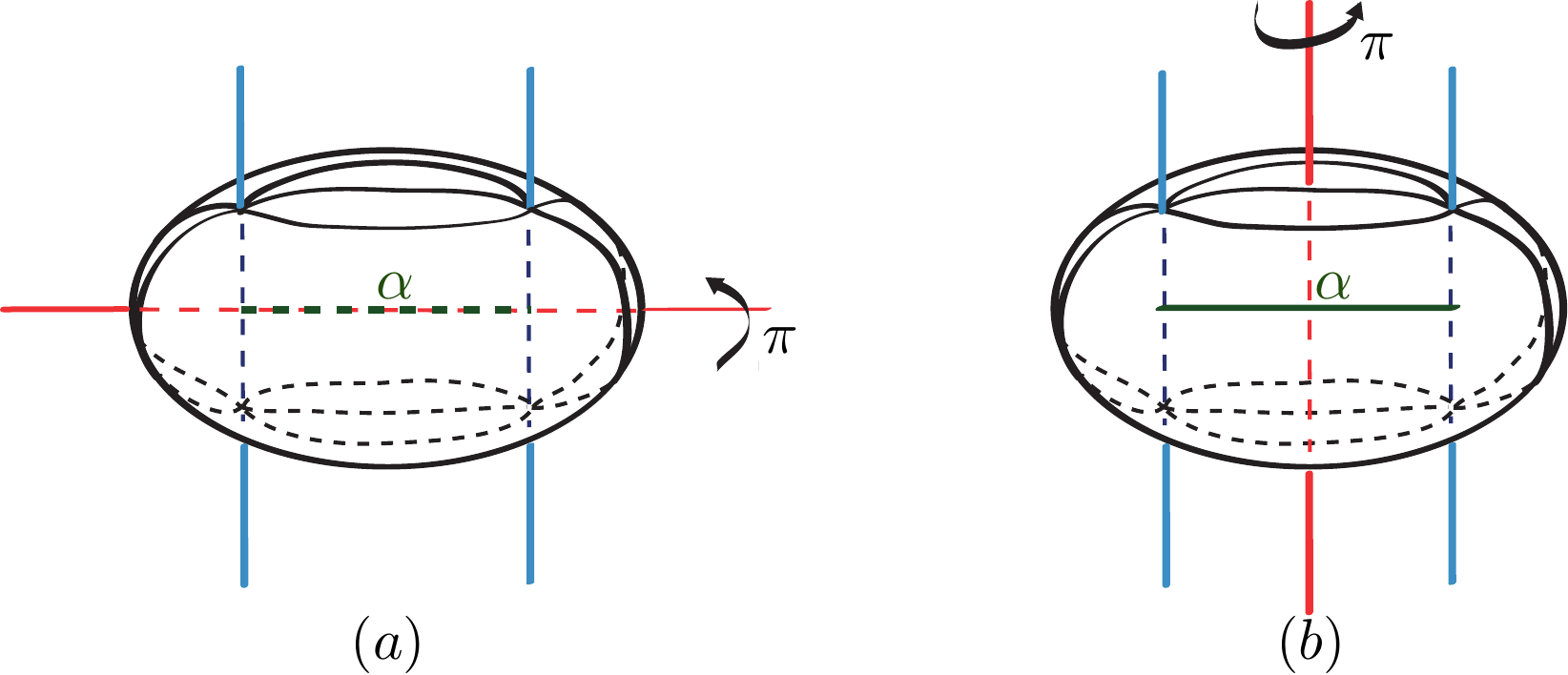}
\caption{Neighborhoods of a proper embedded invariant arc on $(S,c)$.}
       \label{topmodelleri}
      \end{center}
\end{figure}

First note that $\tilde{f}\circ \tilde{c}=\tilde{c}\circ \tilde{f}^{-1}$ and $\tilde{c} \circ \sigma \circ \tilde{c}=\sigma^{-1}$ since $\tilde{c}$ is orientation reversing. 
Then it is easy to see that $(S',\tilde{f}\circ\sigma, \tilde{c}\circ \sigma)$ 
is really an abstract real open book. In fact,
$$
(\tilde{f}\circ \sigma)(\tilde{c}\circ \sigma)=\tilde{f}\circ \tilde{c}\circ (\tilde{c}\circ \sigma\circ  \tilde{c})\circ \sigma = \tilde{c}\circ \tilde{f}^{-1}\circ \sigma^{-1}\circ \sigma = (\tilde{c}\circ \sigma)\circ (\sigma^{-1}\circ \tilde{f}^{-1}). 
$$
Furthermore the real manifold $M'=M_{\tilde{f}\circ \sigma}$ is equivariantly diffeomorphic to $M_f$.
Let us see this in the case where $H$ is a single 1-handle (the other case is similarly treated). Observe that a sufficiently small closed neighborhood $D'$ of the core $\alpha$ of $H$
in $M'$ is a real 3-ball, which is known to be unique up to equivariant  isotopy.
The local model
for $D'$ can be taken as depicted in Figure~\ref{topmodelleri} (see e.g \cite{gg}, Section~1). If $\tilde{c}$ fixes $\alpha$ pointwise
then the model is as in Figure~\ref{topmodelleri}(a) or else $\tilde{c}$ acts on $\alpha$ as reflection and we have the
model as in Figure~\ref{topmodelleri}(b).
A positive stabilization in the non-abstract sense is nothing but excising a sufficiently small closed neighborhood $D$
of $\alpha$ and gluing back a model for $D'$ depending on the behavior of the real structure on $\alpha$.
Thus the uniqueness of the model in each case guarantees that 
the (equivariant) identity  diffeomorphism from $M_f\setminus D$ to $M'\setminus D'$ uniquely extends (equivariantly)
to a diffeomorphism from $D$ to $D'$.

\begin{figure}[h]
   \begin{center}
\includegraphics[scale=0.46]{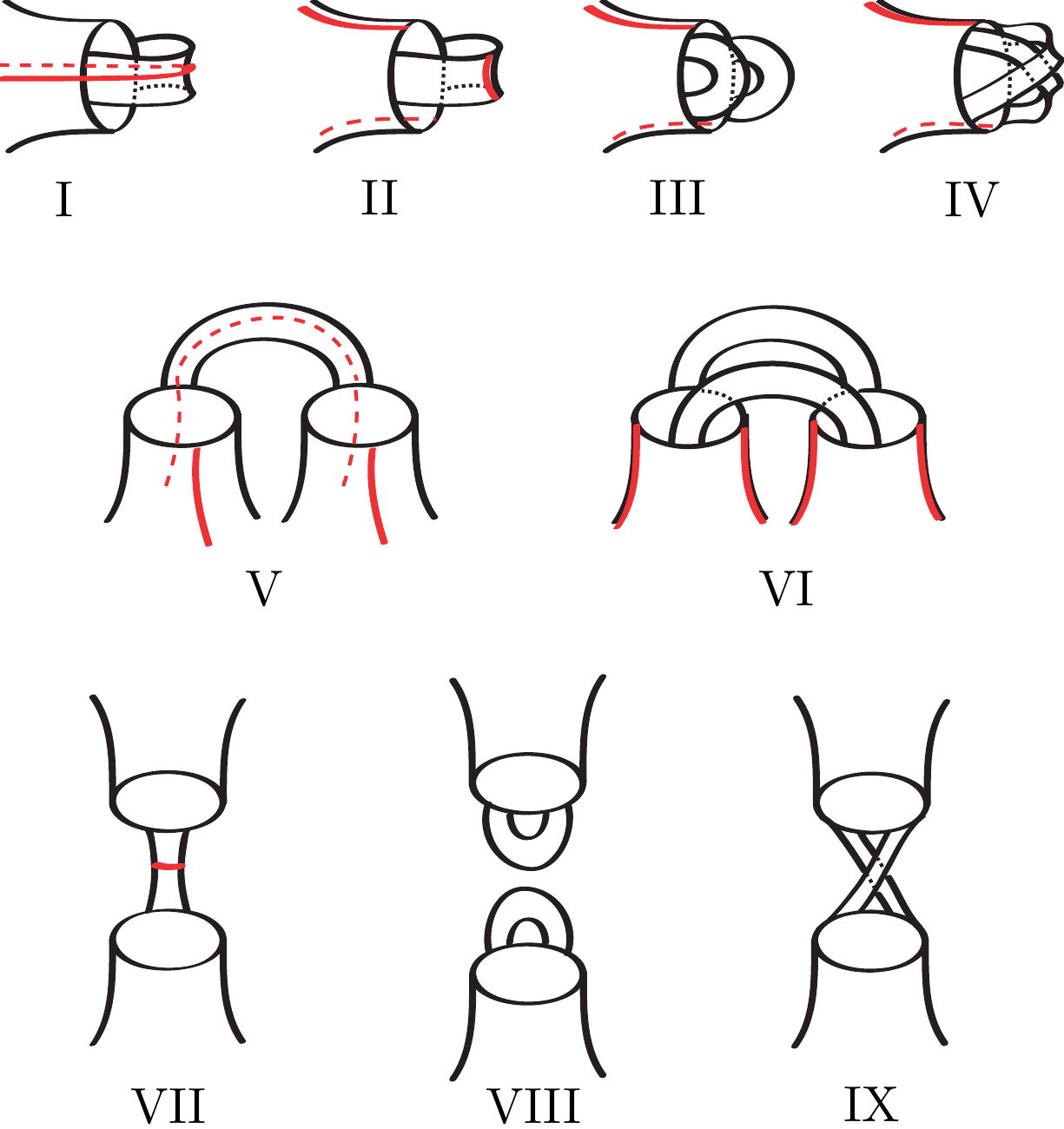}
\caption{Possible handle attachments and extension of the real structure through the handle(s). The real parts
are depicted in red.}
       \label{stablar}
      \end{center}
\end{figure}

In Figure~\ref{stablar}, we depict all possible choices for the attachment(s) of the 1-handle(s) and the extension { of} the real structure through the handle(s). In the first six cases, the { initial} real structure acts as a reflection on the { initial}  boundary component(s) of $S$, while in the last three cases, the { initial}  real structure swaps two boundary components.

\begin{figure}[h]
   \begin{center}
\includegraphics[scale=0.5]{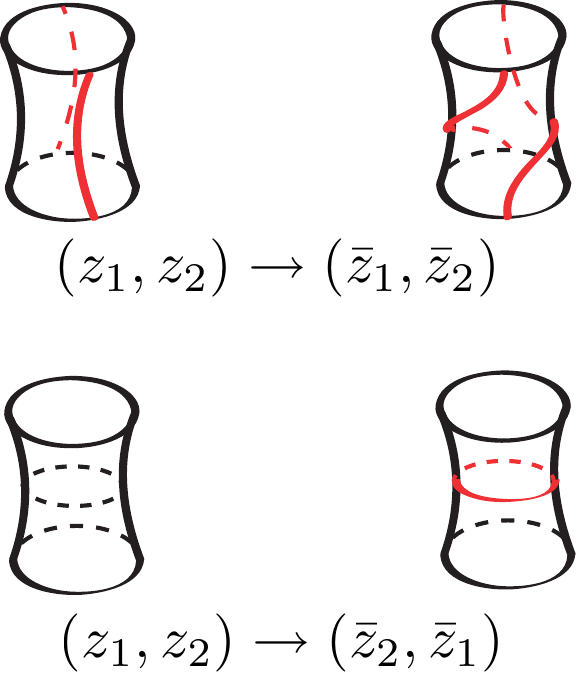}
\caption{Real pages of possible real open books of $S^3$ whose binding is the positive Hopf link.}
       \label{hopfstablari}
      \end{center}
      \end{figure}

\begin{remark} 
A positive real stabilization of type I, II, V, VII is a connected sum with real $S^3$ equipped with the real open book (defined by $\pi: S^3\setminus H\to S^1$, $\pi(z_1,z_2)=\frac{z_1z_2}{|z_1z_2|}$)  whose binding is the positive Hopf link $H=\{(z_{1}, z_{2})\in S^3: z_{1}z_{2}=0\}$, where $S^3=\{(z_{1}, z_{2})\in \C^2: |z_{1}|^2+|z_{1}|^2=1\}$.
This open book 
admits two real structures: $(z_{1}, z_{2})\mapsto (\bar{z}_{1}, \bar{z}_{2})$ and $(z_{1}, z_{2})\mapsto(\bar{z}_{2}, \bar{z}_{1})$.
In the case of $(z_{1}, z_{2})\mapsto (\bar{z}_{1}, \bar{z}_{2})$, the real structure acts on each binding component and the real pages are as depicted at the top of Figure~\ref{hopfstablari}. In the other case, however, the real structure takes one binding circle to the other and acts on the real pages as shown  at the bottom of Figure~\ref{hopfstablari}. 

In \cite{os}, we showed that there is a unique tight contact 3-ball. Using that and the results in the proof of that, it can be deduced that, as in the non-real case, if there exists  a real contact structure   \emph{supported} by a real open book, a positive real stabilization 
of the open book  coincides with the real contact connected sum of the manifold with real $S^3$.

Other types (III, IV, VI, VIII, IX) of positive real stabilizations can be seen as  connected sums with a pair of $S^3$ exchanged by the real structure.
\end{remark}

\subsection{Real Heegaard decompositions}

Here we present the notion of a real Heegaard decomposition which is closely related to a real open book.

Let $\Sigma$ be an oriented hypersurface in a closed oriented manifold $M$. If the topological closure of $M\setminus \Sigma$ is a
disjoint union of two handlebodies diffeomorphic to each other, then the pair $(M,\Sigma)$ is called an \emph{embedded Heegaard decomposition}.

Let $N$ be an oriented handlebody with nonempty boundary $H$ and $h:H\rightarrow H$ be an orientation reversing
diffeomorphism. The pair 
$(H, h)$ is called an \emph{abstract Heegaard decomposition} of the manifold $(N \cup N')/_{x\sim h(x), x\in H}$
where $N'$ is a copy of $N$.

It is obvious that an abstract Heegaard decomposition determines an embedded one immediately.
For the converse, 
consider an embedded Heegaard decomposition  $(M, \Sigma)$, where $\Sigma$ separates $M$ into two diffeomorphic 
handlebodies $N_1$ and $N_2$, lying in $M$.
Furthermore we fix an orientation preserving diffeomorphism $f:M\ra M$ 
such that $f(N_2)=N_1$ and $f|_{\Sigma}$ is orientation reversing.  
Then the pair $(\Sigma, f|_{\Sigma})$ defines an abstract Heegaard decomposition for $M$. 
In fact, suppose $N$ and $N'$ are two identical copies of $N_1$. Then it is straightforward to see that
the map $\varphi_f:M\ra (N\cup N')/_{x\sim f(x),\mbox{ } x\in\Sigma}$ 
defined on $N_1$ as $id:N_1\ra N$ and on $N_2$ as $f|_{N_2}:N_2\ra N'$
is well-defined which is furthermore a diffeomorphism. Observe that here we do not require that $f$ is an involution.

Let $(M,c_{M})$ be a real manifold and $(M, \Sigma)$ be an embedded Heegaard decomposition. If $c_{M}(\Sigma)=\Sigma$ and $c_{M}|_{\Sigma}$ is orientation reversing, then we call the triple $(M, c_{M}, \Sigma)$ an \emph{embedded real Heegaard decomposition}. 
Similarly consider the abstract Heegaard decomposition  $(\Sigma, f)$ with $\Sigma$  bounding the handlebody $N$. 
Consider an orientation preserving involution $c:N\ra N$ such that on $\Sigma$ we have $c\circ f\circ c=f^{-1}$. 
Then we call the triple $(\Sigma,f,c)$ an \emph{abstract real Heegaard decomposition}. 
It is again immediate to see that the embedded real Heegaard decomposition $(M, c_{M}, \Sigma)$ is equivariantly
diffeomorphic to the abstract real Heegaard decomposition  $(\Sigma,c_M|{\Sigma},id)$.
In fact, the map $\varphi_c$ as we have constructed  in the last
paragraph is an equivariant diffeomorphism between these two manifolds.
Finally, note that it was proven by T.~A.~Nagase that every real manifold admits an embedded real Heegaard decomposition \cite{nag}.

Recall that an embedded (respectively abstract) open book decomposition describes an embedded  (respectively abstract) Heegaard decomposition. 
Let $(B,\pi)$ be an embedded real open book for the real manifold $(M,c_{M})$.
Denote by $c_-$ and $c_{+}$ the diffeomorphisms obtained by restricting $c_{M}$ to the two invariant pages $S_-$ and $S_{+}$.
Note that $f=c_{+}\circ c_-$.
Then $(M, S_-\cup_{\partial} S_{+})$ is an embedded real Heegaard decomposition  and 
 $(S_-\cup_{\partial} S_{+},  c_-\cup f\circ c_{-},id)$,
is an abstract real Heegaard decomposition for $(M,c_{M})$.

\section{Real open books in three dimensions}
\label{rak}
\subsection{Real open books on $S^3$}
Here we give a series of examples of real open book decompositions and real  (embedded)  Heegaard decompositions on $S^3$.
It is well-known that on $S^3$, up to isotopy, there is a unique real structure with nonempty real part; its fixed point set is an unknot 
\cite{w} (cf. \cite{bbm}).
Let us fix the real structure $c_{0}$ on $S^3$, considered as one point compactification of $\R^3$,  as the one induced from the rotation on $\R^3$ by $\pi$ along the $y$-axis. We consider the open book on $S^3$ with binding  $z$-axis $\cup\,\infty$ and with pages topologically disks.
Abstractly, this is the open book of $S^3$ with page $S_0$ a disk and monodromy $f_0$ the identity map. 
There are two invariant pages on which the real structure acts as the reflection, $\rho_{0}$, with respect to the $y$-axis.  The corresponding Heegaard decomposition has the splitting surface $\Sigma_0$ a sphere and the gluing map  $\varphi_0$ between the two handlebodies is a reflection on $\Sigma_0$ fixing an \emph{equator}.

Let us perform a positive real stabilization on the real open book we have. At this point it is important to note that we will describe the stabilization on an abstract real open book. However,  one should keep in mind that abstract positive real stabilization can be seen as an \emph {equivariant} connected sum with $(S^3, c_{c_{S^3}})$ equipped with a positive Hopf real open book.

Now, consider the open book $(S_0, f_0)$, we make a positive stabilization of type I (see Figure~\ref{stablar}) to obtain
the page $S_1$ an annulus. The monodromy $f_1$ is the Dehn twist $\tau_1$ along the core $a_1$ of $S_1$.  
The real structure $c_{1-}$ on { $S_{1}$} is the reflection with respect to the core $a_{1}$ composed by the Dehn twist along $a_{1}$, while the real structure $c_{1+}$ on the opposite real page  acts as the reflection with respect to  $a_{1}$;
the corresponding Heegaard surface $\Sigma_1$ is a torus; the real structure on $\Sigma_1$ is $c_1 = c_{1-}\cup c_{1+}$ which 
is exactly the gluing map $\varphi_1$ between the two handlebodies. Note that $c_{1}$  is conjugate to the involution which interchanges the meridian and  the longitude. 

\begin{figure}[h]
   \begin{center}
   \includegraphics[scale=0.7]{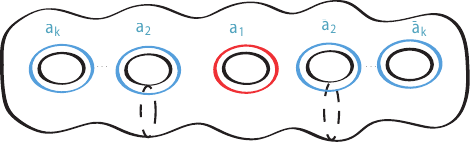}
\caption{Heegaard splitting of odd genus with nonseparating real part.}
       \label{stabor1}
      \end{center}
\end{figure}

In the next steps we perform several positive real stabilizations of type VIII so that 
we obtain a sequence of real open books and real Heegaard decompositions of odd genus for $S^3$: given odd $g=2k-1$,
there is a real Heegaard splitting, see Figure~\ref{stabor1}, 
with surface $\Sigma_g$ of genus $g$ and with the real structure 
and the gluing map given by $\varphi_{g}=c_{k-}\cup c_{k+}$ on $\Sigma_{g}$
where $c_{k-}=\rho_k\circ \tau_1 \circ  \tau_2\circ\tau_{\bar{2}} \ldots  \circ \tau_{\bar{k}}$ with $\rho_{k}$, the real structure on $\Sigma_{g-}$ obtained by extending $\rho_{0}$
to the added handles. 

\begin{figure}[h]
  \begin{center}
\includegraphics[scale=0.7]{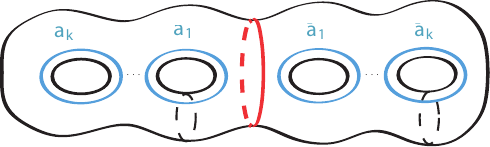}
\caption{Heegaard splitting of even genus with a separating real part.}
       \label{stabor2}
      \end{center}
\end{figure}

Now starting from the real open book $(S_{0}, f_{0}, \rho_{0})$ of $S^3$ with disk pages, we perform a sequence of positive real stabilizations of type III. Thus we get a real Heegaard splitting, see  Figure~\ref{stabor2}, 
with surface $\Sigma_g$ of genus $g=2k$ and with the real structure 
and the gluing map $\varphi_g=c_{k-}\cup c_{k+}$
on $\Sigma_g$ where $c_{k-}=\rho_{k}\circ \tau_1 \circ  \tau_{\bar{1}}\circ  \ldots  \circ \tau_{\bar{k}}$ is the real structure on  $\Sigma_{-}$. Note that in this case the real part of the real structure $\varphi_{g}$ is a separating curve.

\begin{figure}[h]
   \begin{center}
\includegraphics[scale=0.7]{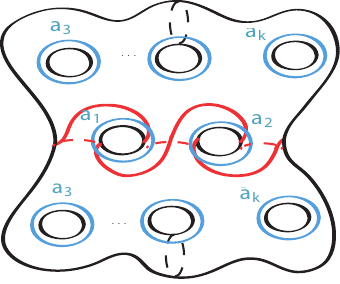}
\caption{Heegaard splitting of even genus with a nonseparating real part.}
       \label{stabor3}
      \end{center}
\end{figure}

Finally, starting from the open book of $S^3$ with disk pages again, we perform twice positive real stabilization of type II and 
then several positive real stabilizations of type III.
We obtain a real Heegaard splitting, see Figure~\ref{stabor3}, 
with surface $\Sigma_{g}$ of even genus $g=2k$ 
and the gluing map $\varphi_g=c_{k-}\cup c_{k+}$, where $c_{k-}=\rho_{k}\circ \tau_1 \circ  \tau_{2}\circ\tau_{3}\circ \tau_{\bar{3}}  \ldots  \circ \tau_{\bar{k}}$.

\subsection{Examples of real Heegaard decompositions on lens spaces}

\label{mercek}

We want to  make some observations on real Heegaard decompositions on lens spaces. 
The classification of real structures on lens spaces up to isotopy is known \cite{hr}.
In op.cit. the real structures of type~$C$ and $C'$ on a lens space are already in the form of real Heegaard decompositions
of genus 1. However there are lens spaces  where there are no real structures of type~$C$ and $C'$  at all 
(see the last three lines of the table in \cite{hr}, Section~4.9). Then the existing real structures on those lens spaces 
cannot be presented as a real Heegaard decomposition of genus 1; however
we know that there {\em is} a real Heegaard decomposition on those lens spaces \cite{nag}.
 Of course, in that case the higher genera decompositions will have the real structure act on the Heegaard surface not maximally. 
We will call a real Heegaard decomposition {\it maximal}
 if the real part of the real structure has genus$+1$ connected components (equivalently, if the number of real components is maximal in view of the Harnack's curve theorem).
This leads to the following observation.
\begin{prop} 
\label{maksimal} 
There are infinitely many real closed oriented 3-manifolds which cannot have a maximal real Heegaard decomposition.
\end{prop}

Several more observations are in order:

\begin{itemize}
\item All open book decompositions with annuli pages can be made real (since the monodromy can only be a power of the Dehn twist 
along the belt circle and any power of the Dehn twist can be written as a product of two real structures). 
Thus the corresponding real manifolds admit real Heegaard splittings of genus 1.
These manifolds fall in two classes: the real lens spaces $L(p,p-1)$  (i.e. the real tight link of the $A_n$ singularities)  and
the real $L(p,1)$ (the supported contact structure being overtwisted). 

\item  The lens spaces $L(p,q)$ where the coefficients of the continued fraction expansion of $-p/q$ are all $-2$ except one which is $-3$
have an open book decomposition with page thrice punctured sphere and with monodromy a product of certain powers of the Dehn twists 
around the boundary components. Since with respect to any real structure on the page the boundary components are invariant,
these open books are always real for any choice of real structure.   Hence, this type of lens spaces have a real genus-2 Heegaard splitting.

\item For an example to above consider $L(5,3)$ on which there are 3 real structures; they are of types~$A$, $B$ and $B'$. 
Hence there is no genus-1 real Heegaard decomposition of $L(5,3)$. Meanwhile, on a sphere with three punctures there are only two real structures. 
Thus, with respect to at least one real structure the corresponding Heegaard decomposition  must have genus greater than $2$.
\end{itemize}

Following this discussion, we ask a natural question:

\noindent {\bf Question.} Given a real closed 3-manifold $(M,c_{M})$, what is the {\em $c_{M}$-real Heegaard genus}, i.e. 
what is the minimum genus among all real Heegaard decompositions representing $c_{M}$?

To finish this section, let us display the question we raised in the Introduction. 

\noindent {\bf Question.} Does every real 3-manifold $(M,c_{M})$ admit a $c_{M}$-real open book?

Note that if the answer is affirmative then with Proposition~\ref{reeltw} below, this would imply that every real 3-manifold  admits a real contact structure.

\section{Real Giroux correspondence}

In this section we give the proof of some parts of
\begin{conj} ({\it Real Giroux Correspondence})
\label{jiru}
Let $(M,c_{M})$ be a real 3-manifold. Then there is a one to one correspondence between the real contact structures on $(M,c_{M})$ up to equivariant contact isotopy
and the real open books on $(M,c_{M})$ up to positive real stabilization and equivariant isotopy.
\end{conj}

We work out the proofs of parts of this conjecture in Propositions~\ref{temasvar}, \ref{aynitemas}, \ref{tasir}. If one can prove Conjecture~\ref{aynikitap}, that will finish the proof of the Real Giroux Correspondence. 

First we explain how to construct a real contact form on the real manifold $(M_{f}, c_{M_{f}})$ built-up as before from the abstract real open book decomposition $(S, f, c)$. 
The construction is the real version of the one of W.~P.~Thurston and H.~Winkelnkemper \cite{tw} and  it 
motivates { the definition of the support relation: 
We say that a real contact structure $\xi=\textrm{ker } \alpha$ on a real 3-manifold $(M, c_{M})$ is {\em supported} by a real open book decomposition $(B, \pi)$ of $(M, c_{M})$ if $\xi$ is supported by the open book $(B, \pi)$  in the usual sense.  Namely, $\alpha$ is positive on the binding considered with its orientation induced from the pages and $d\alpha$ is a symplectic form on each page of the open book.} 

\begin{prop}\label{reeltw}
Every real open book decomposition $(S, f, c)$ supports a real contact structure $\xi=\textrm{ker }\alpha$ on $ M_{f}=S_{f} \bigcup (\amalg_{| \partial S |} S^1\times D^2). $
\label{temasvar}
\end{prop}

\noindent {\it Proof:} The existence proof of the proposition with no equivariance requirement applies with slight modifications to equivariant case. 

Let $\beta$ be a 1-form on $S$ such that $\beta=e^s d\theta$ near $\partial S$ and $d\beta$ is a positive volume form on $S$.  It is elementary to show that such a form exists (see e.g. \cite{g}, Section~4.4.2). Let $\mu(t):[0,1] \to [0,1]$ be a smooth function which is identically 0 near $t=0$, and identically 1 near $t=1$. We consider  the 1-forms 
$$
\begin{array}{l}
\beta^-_{K}=\beta +Kdt \mbox{ on } S \times I_{-}\\ 
\beta^+_{K}=(1-\mu(t))c^*\beta+\mu(t) (f\circ c)^*\beta-Kdt  \mbox{ on }  S\times I_{+}.
\end{array}
$$

Define $\hat{\beta}_{K}$ such that   
$\hat{\beta}_{K}|_{S\times I_+}=\beta^+_{K}$ and 
$\hat{\beta}_{K}|_{S\times I_-}=\beta^-_{K}$ which, by definition, induces a 1-form on the mapping torus  
$$S_f=\big((S\times I_+)\cup(S\times I_-)\big)/_{(x,0)\sim (c(x),0), (x,1)\sim(f\circ c(x),1).}$$ 
However, this form is in general not real, therefore we consider the 1-form $\alpha_{K}=\hat{\beta}_{K} -c_{S_{f}}^*(\hat{\beta}_{K})$. By definition $\alpha_{K}$ is a real form.  Since  $c_{S_{f}}$ is nothing but the  identity map, we have
$\alpha_{K}|_{S\times I_+}= \beta^+_{K}-\beta^-_{K}$, while $\alpha_{K}|_{S\times I_-}= \beta^-_{K}-\beta^+_{K}$.

We claim that for large values of $K$,  $\alpha_{K}$ is a real contact form on $S_{f}$.
To prove this we  show that
$(\beta^+_{K}-\beta^-_{K})\wedge d(\beta^+_{K}-\beta^-_{K})$ is positive on $S_{f}$ for sufficiently large $K$ where positivity is determined with respect to the orientation defined by the form $d\beta \wedge dt$.  

We have
$$
\begin{array}{rcl}
(\beta^+_{K}-\beta^-_{K})\wedge d(\beta^+_{K}-\beta^-_{K}) & = &
2Kdt\wedge d\beta \\
& & -2Kdt\wedge [(1-\mu(t))dc^*\beta+\mu(t) d(f \circ c)^*\beta]\\
& & + \mbox{\footnotesize{ terms not containing }} K \\
& = & 2Kd\beta\wedge dt \\
& & -2K [(1-\mu(t)) c^*d\beta+\mu(t) (f \circ c)^*d\beta]\wedge dt\\
& & + \mbox{\footnotesize{ terms not containing }} K
\end{array}
$$
As $c$ is orientation reversing, both $-c^*d\beta$ and $-(f\circ c)^*d\beta$ are area forms on $S$ 
defining the same orientation as $d\beta$.  Thus,  the first two terms 
above are positive on $S_{f}$ and is dominant when $K$ is large enough.

For the extension of $\alpha_{K}$ over the solid tori, we need to consider two separate cases which are distinguished by  the action of $c_{S}$ on $\partial S$. 
As discussed above, either  $c_{S}$ acts on a boundary component as reflection or it switches two boundary components.

Case 1: Let $S^1_{\partial}$ denote the boundary component of $\partial S$ on which $c_{S}$ acts as reflection.  

We now consider $\nu(S^1_{\partial})=\{(s, \theta): s\in [-\epsilon, 0], \theta \in [-\pi, \pi] \}$ where $c_{S}|_{\nu(S^1_{\partial})}(s, \theta)=(s, -\theta)$.   
On the other hand, we take $S^1\times D^2=\{(\vartheta, r, \varphi): \vartheta, \varphi \in[-\pi, \pi], r\in [0,1] \}$  together with the  real structure $c_{S^1\times D^2}(\vartheta, r, \varphi)=(-\vartheta, r, -\varphi)$ and the 1-form $\alpha'_{r}=h_{1}(r)d\vartheta+h_{2}(r)d\varphi$. 
Note that, $\alpha'_{r}$ is a real form for all $r$ and is a contact form for those $r$ satisfying $h_{1}h'_{2}-h'_{1}h_{2}>0$. 
 
We identify  $S^1\times[1-\epsilon, 1]\times S^1 \subset S^1\times D^2$ with $\nu(S^1_{\partial})\times S^1$  by an equivariant orientation preserving  diffeomorphism 
$\Upsilon$ defined as
\begin{displaymath}
\begin{array}{rcl} 
\Upsilon:S^1\times[1-\epsilon, 1]\times S^1 & \to& \nu(S^1_{\partial}) \times S^1\\
(\vartheta, r, \varphi)&\mapsto&(1-r-\epsilon, \vartheta, \varphi). 
\end{array}
\end{displaymath}
Thus, we get 
\begin{displaymath}
\Upsilon^*(\alpha_{K}|_{\nu(S^1_{\partial}) \times S^1})=\left \{\begin{array}{ll}
-(2e^{1-r-\epsilon} d\vartheta+2Kd\varphi) &  \textrm{if $\varphi  \in [0,\pi], $} \\
2e^{1-r-\epsilon} d\vartheta+2Kd\varphi&  \textrm{if  $\varphi  \in [-\pi,0] .$}
\end{array}\right.
\end{displaymath}

Since we require the extended 1-form to be positive on the binding, and to match with $\alpha_{K}$ on $S^1\times[1-\epsilon, 1]\times S^1$,   it is enough to find smooth $h_{1}$ and $h_{2}$ such that 
\begin{displaymath}
\begin{array}{ll}
(1)\;  h_{1}(r)=1  \textrm{ and } h_{2}(r)=r^2 & \textrm{near }  r=0\\
(2)\;  h_{1}(r)=2e^{1-r-\epsilon}  \textrm{ and }  h_{2}(r)=2K & \textrm{near }  r=1\\
(3)\;  h_{1}(r)h'_{2}(r)-h'_{1}(r)h_{2}(r)>0 & \forall r \in [0, 1],
\end{array}
\end{displaymath}
and it is easy to see that there exist $h_{1}, h_{2}$ satisfying the above conditions.

Case 2:    Let $S^1_{\partial_{1}}$ and $S^1_{\partial_{2}}$  be two boundary components of $S$ such that $c_{S}(S^1_{\partial_{1}})=S^1_{\partial_{2}}$. 
We consider two contact solid tori $S^1\times D^2$ with contact structures $\alpha'(r)={ h_{1}(r)}d\vartheta+ h_{2}(r)d\varphi$ and  $c^*_{S^1\times D^2} (\alpha'(r))=-h_{1}(r)d\vartheta- h_{2}(r)d\varphi$. 

As  before, we identify $S^1\times [1-\epsilon, 1] \times S^1 \subset S^1\times D^2 $ with  $\nu(S^1_{\partial_{1}}) \times S^1$
by $\Upsilon$ then the identification $\Upsilon'$ of  $\nu(S^1_{\partial_{2}}) \times S^1= c_{S_{f}} (\nu(S^1_{\partial_{1}})) \times S^1$ is determined by the symmetry
 \begin{displaymath}
\xymatrix{S^1\times [1-\epsilon, 1]\times S^1  \ar[r]^{\Upsilon} \ar[d]_{c_{S^1\times D^2}} & \nu(S^1_{\partial_{1}})\times S^1 \ar[d]^{c_{S_{f}}} \\
        S^1\times [1-\epsilon, 1]\times S^1   \ar[r]^{\Upsilon'} &\nu(S^1_{\partial_{2}})\times S^1 .}
\end{displaymath}

Near $S^1_{\partial_{1}}$, the extension can be done as in the non-real case 
and the extension of $\alpha_{K}$ near $S^1_{\partial_{1}}$ is obtained symmetrically, using $-h_{1}$ and $-h_{2}$.
\hfill  $\Box$ \\

\begin{prop}
\label{aynitemas}
Two real contact structures supported by the same real open book decomposition are equivariantly isotopic.
\end{prop}

\noindent {\it Proof:} 
Let $ \xi_{0}=\textrm{ker }\alpha_{0}$ and $\xi_{1}=\textrm{ker }\alpha_{1}$ be two real contact structures supported  by the same real open book decomposition. We first isotope $\alpha_{i}, i\in \{0,1\}$ near the binding. 
 
Let $S^1 \times D_{\epsilon}^2=\{(\vartheta, r, \varphi): r\in [0, \epsilon]; \varphi, \vartheta \in [-\pi, \pi]\}$ be  a neighborhood of a binding component, where the real structure is given by $c:(\vartheta,r,\varphi)\mapsto ({-}\vartheta,r,{-}\varphi)$. 
Take a function  $h(r)$  such that $h(r)= r^2$ near 0, $h(r)=1$ near $r=\epsilon$ and $h'(r)>0$, for all $r$.  Define $\alpha_{i,R}=\alpha_{i}+Rh(r)d\varphi$, with  $R\geq0$. 
It can be  checked easily that both $\alpha_{0,R}$ and  $\alpha_{1,R}$ are contact and $c$-real,  and they are equivariantly isotopic to $\alpha_{0}$ and $\alpha_{1}$ respectively. Now,  we consider $\alpha_{t,R}=(1-t)\alpha_{0,R}+ t\alpha_{1,R}$. For each $t\in[0,1]$  this is a real form and is a contact  form   for sufficiently large $R$. 

{ Otherwise if $c$ swaps two distinct binding components then the above modification near one and the $c$-equivariant   modification around the other give the desired isotopy.}   Hence we get an equivariant isotopy between $\alpha_{0}$ and $\alpha_{1}$ on a neighborhood of each binding component.

Furthermore, { the isotopy extends} over all $M$ and hence is an isotopy away from the binding,
between   $\alpha_{0}+R d\varphi$ and $\alpha_{1}+R d\varphi$, which in turn gives an isotopy between $\alpha_0$ and  $\alpha_1$.
Since $c_M$ acts as reflection on $S^1$ (hence sends $\varphi$ to $-\varphi$ outside the binding)  this isotopy is $c_M$-equivariant.
\hfill  $\Box$ \\

\begin{prop}
\label{tasir}
Every real contact structure on a closed 3-dimensional real manifold is supported by a real open book.
\end{prop}

The proof follows from Lemma~\ref{hucreleme} and Lemma~\ref{kurdele aksi}.
First we construct a {\it real contact cell decomposition} on  $M$ (Lemma~\ref{hucreleme})
and then using that we build the real open book (Lemma~\ref{kurdele aksi}).

Let $(M,\xi,c_{M})$ be the real contact manifold. A  {\it real contact cell decomposition} over $M$ is
a  cell decomposition of $M$ with each $k$-skeleton $c_{M}$-invariant ($k=0,1,2,3$);  
each 2-cell is convex with tw$(\partial D,D)=-1$;
and $\xi$ is tight when restricted to each 3-cell.

\begin{lemma}\label{hucreleme}
Every closed real contact 3-manifold $(M,\xi,c_{M})$ has a real contact cell decomposition.
\end{lemma}

\noindent {\it Proof:} 
 { 
We cover $M$ with a $c_{M}$-invariant, finite collection of Darboux balls. 
We then choose  an equivariant cell decomposition such that each 3-cell lies in the interior of a Darboux ball. Note that every real manifold has an equivariant cell decomposition \cite{il}. 
We turn the equivariant cell decomposition into an equivariant contact one: { we first make the 1-skeleton Legendrian, then we make the faces convex. Finally we make sure that each face has twisting $-1$.}

 We follow the detailed elucidation in \cite{vau}. To start with, we note that every construction there, provided that it is made away from the real points, can be made equivariantly in our case with its symmetric counterpart. Thus we focus on the construction near real points.
 
{ There is no additional difficulty in the real case to make 1-skeleton Legendrian. Thus, Etape.1  in \cite{vau} can be followed in the  equivariant setting as well; all possible cases in the equivariant contact local picture (which is known to be unique by \cite{os}) is depicted in Figure~\ref{lok}.}
Having this local model in mind, we make each 1-cell Legendrian.
The real edges are already Legendrian. If an edge has a single real point in its interior, the edge can be perturbed 
equivariantly to a Legendrian one. As for a neighborhood of a real vertex, it is obvious that an equivariant pair of edges can be made Legendrian symmetrically.
\begin{figure}[h]
   \begin{center}
\resizebox{12cm}{!}
{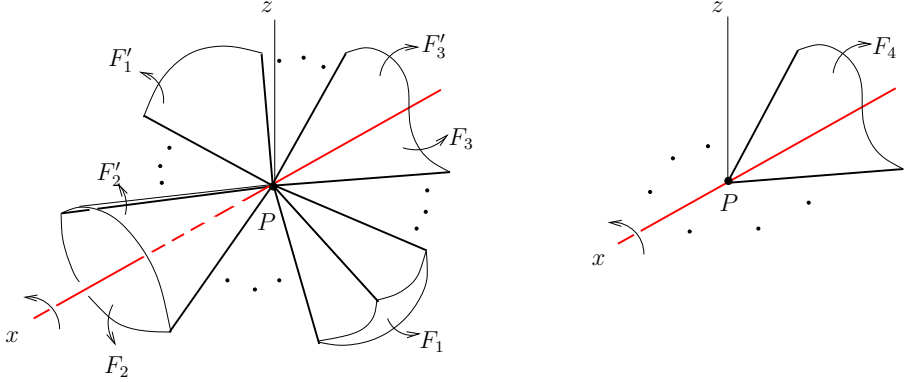}
\caption{All possible sorts of cases for an equivariant cell decomposition near a real vertex $P$: either the 2-cells come in equivariant pairs as $(F_j,F'_j)$, $1\leq j\leq 3$, or contain a real arc as $F_4$. All edges lie on $xy$-plane}
     \label{lok}
      \end{center}
\end{figure}

{
Now we modify the 2-skeleton to get convex 2-cells (in \cite{vau}, Etape.2 and Etape.3 achieve that carefully).  
As each 2-cell remains inside a Darboux ball,  the (weak) Bennequin inequality applies, so  tw$(\partial b, b)\leq -1$ for each 2-cell $e$. Hence we might hope that we can perturb each $e$  convex as in \cite[Proposition~3.1]{ho}. However more care is needed here since the Legendrian boundaries are not smooth. 
Following \cite{vau} we first make the faces well-positioned with respect to $\xi$.
Namely, we make sure that there are no two faces of a 3-cell sharing an edge such that this common edge contains an arc $[b,c]$ with $b$ a negative singularity on one face and $c$ a positive singularity on the other face, and $(b,c)$ is a leaf of the characteristic foliation of the two faces oriented from $b$ to $c$.  
Below we refer to  the faces with no such incidence as {\em $\xi$-good}. }
To apply the procedure in \cite{vau}, 
one should start by getting 
tw$(a)<0$ (relative end points) for every 1-cell $a$. This is possible for every 1-cell in the given 
cell decomposition except possibly the real edges. Note that a real edge cannot be isotoped to 
a different set equivariantly. Therefore if a real edge $r$ does not already have a negative twisting,
$r$ must be deleted from the cell decomposition. One way to do that is as follows: let exactly $k$ faces
have $r$ in their boundary. In each face, take an arc with the same end points as $r$ such that 
the arcs in symmetric faces are chosen symmetric. 
If $k=2$, take an extra, equivariant pair of arcs with the same end points and not
intersecting the interior of  any face and set $k=4$. 
Let  us index those arcs in the order they appear revolving around $r$. Suppose that $r_j$'s are chosen very close to $r$;  that $r_j$'s have no real points in their interior; and that each pair $r_i,r_j$ of  
cyclically  consecutive indices  bounds an embedded disk $e_{ij}$ such that these $k$ disks are disjoint, they have no real points in their interior and they come in equivariant pairs, 
i.e. $c_M(e_{ij})$ is the disk across. 
We delete $r$ and include $r_1,\ldots,r_k$ to the 1-skeleton, $e_{ij}$'s to the 2-skeleton and
the 3-ball $B$ enclosed by $e_{ij}$'s in the 3-skeleton (see Figure~\ref{reelyaysil} for initial $k=2$). 
Moreover the
2- and 3-cells touching $r$ are retracted off $e_{ij}$'s and $B$ accordingly. Note that $B$ 
lies in a Darboux ball provided that $r_j$'s are sufficiently close to $r$.
After this modification in the cell decomposition, we go back to the last paragraph to modify 
the 1- and 2-skeletons around the real vertices to get the local model again.
\begin{figure}[h]
   \begin{center}
\resizebox{6cm}{!}
{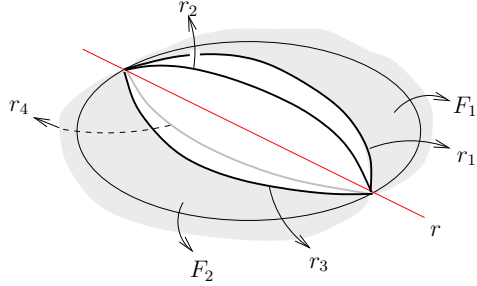}
\caption{Deleting a real edge $r$ that lies in the boundary of two 2-cells $F_1$ and $F_2$.}
     \label{reelyaysil}
      \end{center}
\end{figure}

Having made this observation, one can follow the procedure described in 
Etape.2 in \cite{vau}  to obtain $\xi$-good configuration of faces. That procedure applies  without any change for an equivariant pair of 2-cells near a real vertex.
It also applies for $2k$ equivariant pairs of 2-cells $(F_1,F_k),\ldots,(F_{k-1},F_{2k})$ near 
a real edge $r$, with the following remark:
in the procedure in Etape.2, all the  singularities of all $F_j$'s along $r$ are made isolated.
This cannot be achieved in the equivariant case, the singularities of $F_j$ ($j\leq k$) and $F_{j+k}$ must occur at the same points on $r$. Still, it is easy to observe that the perturbations can be done equivariantly and thus the procedure applies without any change.  Note that the procedure needs strictly negative twisting on $r$.

The remaining two cases that involve real points are when $c_{M}$ acts on a 2-cell $e$ as rotation  or as reflection. {First we argue that 
there can be no equivariant retrograde saddle-saddle connection, i.e. from a negative one to a positive one.
(Note that  the nonequivariant ones are ruled out by the discussion in Etape.2.) }
In the former case, { since the real structure preserves the orientation} 
there is no equivariant saddle-saddle connection on $\partial e$
and thus $e$ can be perturbed symmetrically to satisfy the $\xi$-good condition  along $\partial e$. 
As for the latter case,
we first modify $e$ equivariantly around its real boundary points so that those 
real boundary points become sinks and sources.  
Then $e$ can be perturbed equivariantly to satisfy the $\xi$-good condition. 

After this modification of the 2-skeleton near the 1-skeleton, one can make the 2-cells convex. 
There are three cases as above. The case of an equivariant pair of 2-cells is handled symmetrically (see \cite[Section~4.8]{g} or \cite[Etape.3]{vau}). In the case $c_{M}$ acts on a 2-cell $e$ as reflection, 
if there is no retrograde saddle-saddle connection on the real arc $r$,
one can perturb the two  halves of $e$ symmetrically to make them convex. More explicitly, since $c_M$ is orientation  reversing on $e$, it preserves the orientation of  the characteristic foliation so that a quotient field is induced on  $e/c_M$. Thus $e$ can be made convex.

{
Otherwise let us assume that there is a retrograde saddle-saddle connection on the real arc $r$. 
There is no obvious reason that rules out this case. In fact a local model 
can be constructed as follows. In $\R^3$, consider the contact forms 
$$\alpha_1=\sin(\pi y) dx + 2 \sin(\pi x) dy + (2 \cos(\pi x)  - \cos(\pi y)) dz$$ 
and 
$$\alpha_2=\sin(\pi x) dy + (1 - \cos(\pi y)) \sin(\pi y) dx + (\cos(\pi y)- \cos(2\pi y) - \cos(\pi x)) dz $$
of Example~4.6.15 in \cite{g}. 
Both are real with respect to the real structure $\pi_y$ that is a rotation by  
$\pi$  around the $y$-axis.  The ($\pi_y$-symmetric) characteristic foliation 
$X_i$ ($i=1,2$) determined by $\alpha_i$ on the $xy$-plane has a saddle singularity $S_i$ at the origin, which is a positive (resp. negative) saddle for $i=1$ (resp. $i=2$). Place the local model
for $S_i$ at the point $(0,(-1)^i,0)$. We rectify both vector fields $\pi_y$-symmetrically in
$H=\{(x,y)|-\delta<x<\delta, -\frac{1}{2}<y<\frac{1}{2}\}$ for $\delta$ sufficiently
small so that for some $u>0$ they both equal $(0,u)$ identically over the $x$-axis and 
that div$_{X_i}(dx\wedge dy)=0$ in $H$ on and only on the $x$-axis. 
First we modify $X_1$ over $\{0\}\times[0,\frac{1}{2}]$ as  follows. 
Let $X_1(x,y)=(P(x,y),Q(x,y))$. Note  $X(0,1/2)=(0,v)$ for some $v>0$; $P_x(0,1/2)=0$ and $Q_y(0,1/2)\neq 0$. We choose $u>v$ and extrapolate $Q$ over $\{0\}\times[0,\frac{1}{2}]$ in a $C^2$-smooth way such that $Q_y=0$ if and only if $y=0$. Now  we extend this vector field over $H$ $\pi_y$-equivariantly so that $X_1(x,0)=(0,u)$. One can take $\delta$ sufficiently small such that 
div$_{X_i}(dx\wedge dy)=P_x+Q_y\neq 0$ except on the $x$-axis. After a similar modification for 
$X_2$, we get an equivariant vector field $Y$ on $H$. A careful inspection of Lemma~4.6.3 in \cite{g}
in the real setting shows that $Y$ is the characteristic vector field for some $\pi_y$-real contact structure in a neighborhood of $H$ in $\R^3$. 
In this way, we get the $\pi_y$-real $y$-axis carrying a $\pi_y$-symmetric  retrograde saddle-saddle connection on $H$.
}

Now, we cannot get rid of a  retrograde saddle-saddle connection on the real arc $r$  by an 
equivariant perturbation, unless tw$(r)< 0$ relative to the end points in the case of which 
one can employ the procedure in \cite[Etape.2]{vau} for $k=2$.  
So if  $r$ has a  retrograde saddle-saddle connection with  tw$(r)\geq 0$ then $r$ should 
not lie in any 2-cell $e$. To dismiss $r$, one can employ Figure~\ref{reelyaysil}, introducing 4 new 1-
cells  $r_1,\ldots,r_4$, 4 new 2-cells and 1 new 3-cell. 
Note that for $r_i$'s close to $r$, the new 3-ball remains in a Darboux ball as well.
Moreover, the end points $P,Q$ of $r$ are added to the 0-skeleton, the 3-balls having $e$ on their boundary are retracted and  the region enclosed by $r_1$ and $r_3$ is excised from $e$ to obtain 2 new faces.
After this modification in the cell decomposition, one has to go 
back and perform the appropriate modifications  for the local model around the real vertices $P$ and $Q$, for Legendrian edges and for $\xi$-goodness.

In the last case where $c_{M}$ acts on $e$ as 
rotation, the characteristic field on $e$ is an anti-symmetric vector field tangent to the boundary with finite number of  hyperbolic singularities. 
The work in \cite{bn} shows that among all such vector fields, the structural stable ones are dense.
Thus we conclude that $e$ can be perturbed convex equivariantly, employing similar
final steps as in the proof of   Proposition~3.3 in \cite{os}.

Once every 2-cell is made convex, we subdivide those 2-cells with tw$(\partial e^2, e^2)= -n<-1$ to get 2-cells with correct twisting $-1$. 
The subdivision of an equivariant pair of 2-cells  can be made symmetrically. This is true even if the 2-cell is $c_M$-invariant since the dividing set can be chosen invariant \cite[Theorem~3.1]{os}. By the Legendrian realization
principle \cite{ho} each edge $b$ added in the 1-skeleton during the subdivision can be made Legendrian keeping the 2-cells convex (if $b$ is a real arc, it is already Legendrian). 
If an end point $p$ of $b$ is not in the 0-skeleton, we add  $p$ to the 0-skeleton. 
If $p$ is already in the 0-skeleton then either $p$ is real and in the local picture there has to appear the symmetric copy of $b$ as well or else $p$ is not real and there has to appear the symmetric copy of $b$ at the vertex $c_M(p)$.

Thus, we obtain a $c_{M}$-equivariant contact cell decomposition on $M$. \hfill  $\Box$

Note that the local model of the 2-skeleton near a vertex  can be  
modified further as in \cite[Etape.5]{vau} so that the edges   (denoted $b$) 
which were added in the 1-skeleton in the paragraph above can be considered in the $xy$-plane. 
Hence after introducing these last edges, we may suppose that Figure~\ref{lok} still depicts the most
general case.
}

\begin{figure}[h]
   \begin{center}
\resizebox{12cm}{!}
{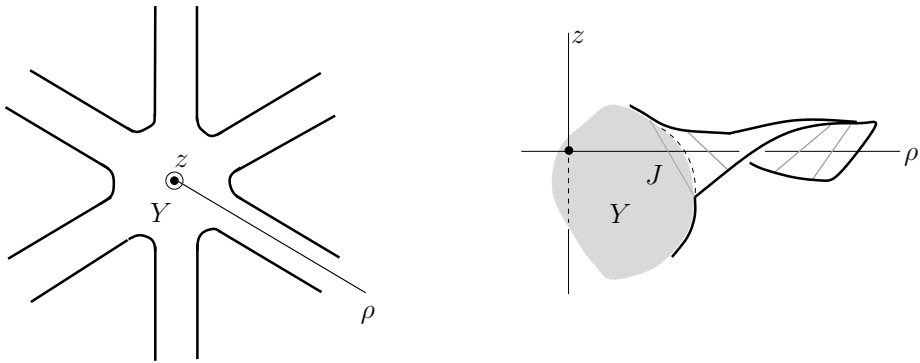}
\caption{The (real) ribbon near a vertex  of valency 6; left: top view, right: the wing along the ray $\rho$.}
     \label{kurdela}
      \end{center}
\end{figure}

Given $(M,\xi,c_{M})$ take a real contact cell decomposition. 
A {\it real ribbon} of the 1-skeleton $G$ is a
compact surface $R$ with boundary such that {
(i) $c_{M}$ acts on $R$ reversing the orientation, (ii) $R$ retracts onto  $G$ 
and (iii) $T_p R$ is arbitrarily close to $\xi_p$ along $G$.

By taking a sufficiently small equivariant regular neighborhood of $G$, one can choose an equivariant, properly embedded surface in this neighborhood which satisfies the conditions above.
Furthermore one can always find a ribbon for which $\xi_p=T_p R$ on $G$ away from the vertices of $G$. In fact, 
the standard neighborhood of a Legendrian arc presently gives a local model for the ribbon away from the vertices of $G$. As for the neighborhood of a { (real)} vertex of $G$ of valency $k$,  we consider the  real Darboux ball near the origin in $\R^3$ with $\pi$-rotation about $x$-axis,
$$(D^3,\beta=dz+r^2d\theta, c:(z,r,\theta)\mapsto (-z,r,-\theta)),$$
and the surface $W \subset D^3$ that is the union of a small disk $Y$ in $xy$-plane centered at 0 and  $k$ radially directed {\em wings} $W_1,\ldots,W_k$ as in Figure~\ref{kurdela}(a). The wings are such that each contains a part of a ray out of 0 in $xy$-plane. For instance one can think of the following  parametrization for $W_j$ and then round the corners. Consider a small secant $J$ of $Y$ with its midpoint on a ray $\rho$ and perpendicular to $\rho$.  The wing $W_j$ is ruled by the family $J_t$ obtained by rotating $J$ around its midpoint faster than $\eta=\ker\beta$ while moving along $\rho$. Once some  $J_{t_0}$ lies in $\eta$, we rotate the line segment with $\eta$. The wing $W_j$ traversed in this way is tangent to $\eta$ along $\rho$ from $t_0$ on and one can check that those are the only  points of $W_j$ where $W$ is tangent to $\eta$. The only point of tangency in  $Y$ is the origin.

In the real case, the wings are to be chosen equivariant altogether.
This is a local model for the (real or not) ribbon near a vertex into which no matter how many edges go. Observe that  $d\alpha=2dx\wedge dy$ restricts to an area form on $W$ and the vector field $\partial_z$ is a contact vector field transverse to $W$, making it convex. The dividing set is empty. Note also that $\partial W$, and  by definition $\partial R$, is transverse to the contact planes. 

Through these observations one can immediately see that the conditions (required in \cite[pages~22,23]{et}) for the dividing set on a ribbon and the relation with the Reeb vector field are satisfied. This latter was expected since we know that the support relation is equivalent to existence of an isotopy having the contact planes arbitrarily close to the tangent planes in compact domains of the pages (see e.g. \cite[Lemma~3.5]{et}).}

{
Note that our model for a ribbon is exactly the one constructed in \cite[Etape.5]{vau}
with the slight modification we have proposed just after the proof of Lemma~\ref{hucreleme}.  
}

{  Now we are ready to state}
\begin{lemma}\label{kurdele aksi}
The boundary $B$ of a real ribbon $R$ in $(M,c_M)$ is the binding of a real open book that supports the real contact structure $\xi$.
\end{lemma}

\noindent {\it Proof:} 
Note that the support relation holds regardless of the real structure and a careful proof has been given in \cite{vau}. 
We assume that we are given the open book supporting $\xi$ such that
the real ribbon is a page, as constructed in \cite{vau} (see also \cite{gi},\cite{et}).
Below we will deform this open book by an isotopy so that it becomes a real open book.

To this end, let us first  choose an invariant neighborhood $N(R)$ of the real ribbon $R$ such that $\partial R\subset \partial N(R)$ and that the boundary $\partial N(R)$ is a convex surface with dividing curve $\partial R$. { A model for such a neighborhood around a vertex of $G$ is
the (smoothed-out) union of small cylinders properly containing wings $W_j$ and a small sphere assuming $Y$ as its equatorial disk. It is easy to see that $\partial Y$ is a dividing set of the sphere. Straightforward gluing techniques of dividing curves reveal that $\partial R$ is a dividing set on  $N(R)$.

Now, w}e can identify $N(R)$ with 
$$\left.\raisebox{.2em}{$R\times [-1,1]$}\middle/\raisebox{-.2em}{$(r, t)\sim (r, t'), r \in \partial R; \; t, t' \in [-1,1]$}\right.$$
on which the real structure acts as $(r,t)\mapsto (c_{M}|_{R}(r), -t)$. Then the projection $R\times (-1,1) /_{\sim}\to (-1, 1)$ defines one ``half" of the open book.  By the construction of the contact cell decomposition each 2-cell $D$ of the cell decomposition is convex and has  boundary with twisting $-1$ relative $D$. Since $R$ twists with the contact structure along the 1-skeleton, each 2-cell $D$ may be assumed to intersect $\partial R$ at exactly two points 
{ (As noted in Footnote~5 in \cite[page~33]{vau}, there might appear additional equivariant pairs of intersection points 
because of the new edges added to the 1-skeleton as in the previous proof. 
In that case one needs to subdivide the 2-cells along equivariant pairs of Legendrian curves of zero tw.).} 
The 2-cells are foliated via the pages of the open book. Each page traces an arc  and all arcs meet at points touching $\partial R$. 

\begin{figure}[h]
   \begin{center}
\includegraphics[scale=0.6]{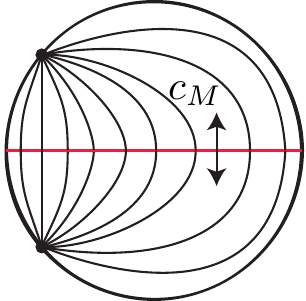}\hspace{1cm}\includegraphics[scale=0.6]{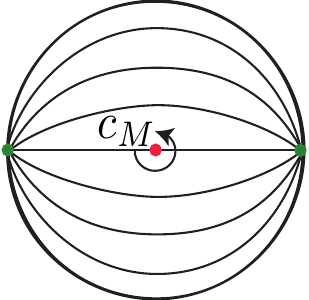}
\caption{A 2-cell foliated by the pages of the open book (left), and possible actions of a real structure on a foliated 2-cell (right).}
     \label{2hucre}
      \end{center}
\end{figure}

\begin{figure}[h]
   \begin{center}
\includegraphics[scale=0.6]{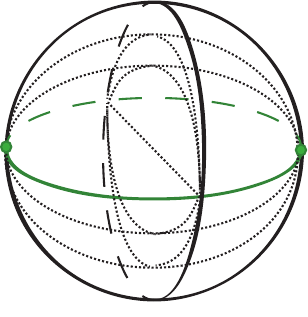}
\caption{A 3-cell foliated by the pages of the open book.}
       \label{3top}
      \end{center}
\end{figure}

Let us now consider the real structure, and first assume that the 2-cell $D$ is invariant, so the boundary of the 2-cell admits an involution which is a reflection or a rotation. In the case of reflection, the involution fixes setwise the points touching $\partial R$. Since any two involutions on a circle fixing the same pair of points can be extended to an involution of a disk in such a way that the extension is unique up to isotopy relative to boundary, we can isotope the open book in a neighborhood of the 2-cell in such a way that it  commutes with the fibration (Figure~\ref{2hucre}, left). In the case of rotation the points touching $\partial R$ are interchanged by the real structure, so  a similar idea applies (Figure~\ref{2hucre}, right). The case when the real structure swaps two 2-cells can also be similarly treated.   
Now, if we trace the dividing curve on the boundary of  each (3-cell$\setminus$int$(N(R))$), we obtain a connected closed invariant curve \cite{gi3}. Therefore,  we can isotope the open book further in the 3-cells in such a way that it
becomes equivariant and the pages foliate the 3-cells as shown in Figure~\ref{3top}. In fact, an involution on $S^2$ with the prescribed image of a fixed closed curve can be extended to a unique (up to isotopy relative boundary) involution of $B^3$ preserving the fibration described above.  
\hfill  $\Box$ \\

To complete the proof of the claimed Real Giroux Correspondence, one needs to prove the following
\begin{conj}
\label{aynikitap}
Two real open books on the real contact manifold $(M,\xi,c_{M})$ supporting the real contact structure $\xi$ 
are related by positive real stabilizations and equivariant isotopy.
\end{conj}

\noindent {\it Idea of the proof:}
The proof will follow the idea of the non-real case proposed by E.~Giroux \cite{gi1}. 
There are two steps in the proof:
showing that any real open book supporting $\xi$ comes from a real contact cell decomposition 
after possibly a number of positive real stabilizations;
and that two such real contact cell decompositions are related by positive real stabilizations.
The latter discussion follows closely the work of L.~Siebenmann \cite{si} and at least the content of
that work can be repeated in the real setting.
\hfill  $\Box$ \\

{
It is worth emphasizing that in the proof of Lemma~\ref{hucreleme}, 
to turn  an equivariant cell decomposition into  a contact cell decomposition, a real edge is
deleted from the 1-skeleton in case it has nonnegative twisting.
This is not unexpected if Conjecture~\ref{aynikitap} holds true. Given a real open book, 
it is not true that there is a series of positive real  stabilizations which produces 
a real open book with the real part lying entirely on a page.}  
For example, the real Hopf open book  depicted at the top of Figure~\ref{hopfstablari} cannot be 
positively equivariantly stabilized so that the real part  becomes a subset of a page. In fact, it is 
immediate to see that if a page contains more than one piece of the same connected component of the 
real  part, then there is no way to get the real part lie on a page by positive real stabilizations.

\end{document}